\documentclass[a4paper,12pt]{article}
\usepackage[utf8]{inputenc}
\usepackage[T1]{fontenc}
\usepackage{amsbsy, amsfonts, amsmath, amssymb, amsthm}
\usepackage{bbm}
\usepackage{calc}
\usepackage{changepage}
\usepackage{enumerate}
\usepackage{etoolbox}
\usepackage{graphicx}
\usepackage{import}
\usepackage{latexsym}
\usepackage{lipsum}
\usepackage{listings}
\usepackage{mathtools}
\usepackage{multicol}
\usepackage{multirow}
\usepackage{setspace}
\usepackage{stmaryrd}
\usepackage{tcolorbox}
\usepackage{tikz}
\usepackage{tikz-cd}
\usepackage{xcolor}
\usepackage{xspace}
\usepackage{hyperref}
\hypersetup{linktoc=all}
\usepackage[noabbrev,capitalize]{cleveref}
\usepackage[export]{adjustbox}
\usepackage{todonotes}
\usepackage{authblk}

\title{ A Generalisation of Niven's Theorem for Trigonometric Functions}

\author{Adam Keilthy \and Ailbhe N{\'i} Ruair{\'i}}

\newcommand{\Addresses}{{
  \bigskip
  \footnotesize

  A.~Keilthy, \textsc{School of Mathematics, Trinity College Dublin, 17 Westland Row, Dublin 2,}\par\nopagebreak \texttt{akeilthy@tcd.ie}

  \medskip

  A. N{\'i} Ruair{\'i}, \textsc{School of Mathematics, Trinity College Dublin, 17 Westland Row, Dublin 2}\par\nopagebreak
  \texttt{rogersa2@tcd.ie}

}}

\date{}

\usepackage[backend=biber,style=numeric-comp,giveninits,url=false,sortcites=true,maxbibnames=99]{biblatex}
\renewbibmacro{in:}{%
    \ifentrytype{article}{}%
    {\printtext{\bibstring{in}\intitlepunct}}}
\def\NN{{\mathbb{N}}}
\def\ZZ{{\mathbb{Z}}}
\def\QQ{{\mathbb{Q}}}
\def\RR{{\mathbb{R}}}
\def\CC{{\mathbb{C}}}
\DeclareMathOperator{\Gal}{Gal}

\def\cosab{{\cos\left(\frac{a\pi}{b}\right)}}
\def\sinab{{\sin\left(\frac{a\pi}{b}\right)}}
\def\tanab{{\tan\left(\frac{a\pi}{b}\right)}}

\def\half{{\frac{1}{2}}}
\def\sqm{{\sqrt[n]{m}}}
\def\tif{{\text{ if } }}
\def\nth{{n^{\text{\normalfont{th}}}}}

\newcounter{mainthm}[section]

\newcounter{mainthrm}[section]

\newtheorem{defn}[mainthm]{Definition}
\newtheorem{prop}[mainthm]{Proposition}
\newtheorem{lem}[mainthm]{Lemma}
\newtheorem{cor}[mainthm]{Corollary}
\newtheorem{thm}[mainthm]{Theorem}

\newtheorem{rem}[mainthm]{Remark}

\definecolor{codegreen}{rgb}{0,0.6,0}
\definecolor{codegray}{rgb}{0.5,0.5,0.5}
\definecolor{backcolour}{rgb}{0.95,0.95,0.92}
\definecolor{codeblue}{rgb}{0.2,0.5,1}
\definecolor{codepurple}{rgb}{0.5,0.2,1}

\lstdefinestyle{mystyle}{
    backgroundcolor=\color{backcolour},  
    commentstyle=\color{codegreen},
    keywordstyle=\color{codepurple},
    numberstyle=\tiny\color{codegray},
    stringstyle=\color{codeblue},
    basicstyle=\ttfamily\footnotesize,
    breakatwhitespace=false,         
    breaklines=true,                 
    captionpos=b,                    
    keepspaces=true,                 
    numbers=left,                    
    numbersep=5pt,                  
    showspaces=false,                
    showstringspaces=false,
    showtabs=false,                  
    tabsize=2  }

\begin{document}

\maketitle
\begin{abstract}
Niven's Theorem asserts that $\{\cos(r\pi)|r\in \QQ\}\cap\QQ = \{0, \pm1, \pm\half\}$. This paper uses elementary methods to classify all elements in the sets $\{\cos^n(r\pi)|r\in \QQ, n \in \NN\}\cap\QQ$ and $\{\sin^n(r\pi)|r\in \QQ, n \in \NN\}\cap\QQ$. Using some algebraic number theory, we extend this to a classification of all elements in $\{\tan^n(r\pi)|r\in \QQ, n \in \NN\}\cap\QQ$. Finally, we present a short Galois theoretic argument to provide a more conceptual understanding of the results. 
\end{abstract}

\pagenumbering{arabic}

\tableofcontents

\newpage

\section{Introduction}

While usually transcendental, the values of trigonometric functions at rational multiples of $\pi$ are easily seen to be algebraic. Niven's Theorem \cite{niven,dresden2009new} shows that the only rational values taken by sine and cosine functions at rational multiples of $\pi$ are the familiar values:
\begin{align*}
    \sin(0) = \cos\left(\frac{\pi}{2}\right) &= 0,\\
    \sin\left(\frac{\pi}{2}\right) =\cos(0) &= 1,\\
    \sin\left(\frac{\pi}{6}\right) = \cos\left(\frac{\pi}{3}\right) &= \frac{1}{2},
\end{align*}
the corresponding negatives, and their shifts by integer multiples of $\pi$. 

The other standard values of these trigonometric functions, namely $\pm\frac{\sqrt{3}}{2}$ and $\pm\frac{\sqrt{2}}{2}$, can be shown to be the only values for which $\cos(r\pi),\,r\in\QQ$, is the square root of some rational number \cite{jahnel}.

Understanding these special values, and the relations between them is both interesting from a purely number theoretic perspective \cite{berger}, and in Euclidean geometry \cite{bergerTriangle}, but has also found applications in hyperbolic Fourier analysis \cite{bylehn}, where the number of radical values taken by trigonometric functions is essential for the validity of certain asymptotic analysis. Algebraic values of trigonometric functions can also be studied through the lens of dynamical systems \cite{panraksa} and transcendence theory \cite{adamczewski2025algebraic}.

A natural question is then to ask whether we can generalise Niven's Theorem to an $\nth$ power. While techniques from Galois theory provide a powerful tool for such an generalisation, we restrict ourselves to more elementary methods. Part of the appeal of Niven's theorem is the simplicity of the statement and the accessibility of the proof, and we aim to maintain this accessibility. 

Thus, we aim to determine whether there are any rational angles $r\pi$ such that $\cos(r\pi)$ is an $\nth$ root of some rational number for $n\geq3,n\in\NN = \{1,2,...\}$. Using results from Berger \cite{berger} relating to the minimal polynomial of $\pm2\cos(r\pi)$, we are able to provide a comprehensive list of all values for which, given $r\in\QQ$ and $n \in\NN$, $\text{cos}^n(r\pi)$ and $\text{sin}^n(r\pi)$ are rational, as outlined in Figure 1 below. Serendipitously, we find that these simple radical values are exactly those familiar to us from introductory trigonometry.

We obtain analogous results for the tangent function, but must appeal to slightly more advanced methods from field theory and basic algebraic number theory, for reasons outlined in Section \ref{sec:tangent} below. Finally, we break from our commitment to elementary discussions to provide a more conceptual explanation for our results using Galois theory, similar to those of \cite{shibukawa2020rational}, for the sake of completeness.
\pagebreak

\renewcommand{\arraystretch}{2}
\begin{center}
    \begin{tabular}{ |c||c|c|c| } 
        \hline
            $b$ & $\cosab$ & $\sinab$ & $\tanab$ \\ 
        \hline
            \vspace{-5.25ex} & & & \\
        \hline
            $1$ & $\pm 1$ & $0$ & $0$ \\
        \hline
            $2$ & $0$ & $\pm1$ & \it undefined \normalfont \\
        \hline
            $3$ & $\pm\half$ & $\pm\frac{\sqrt{3}}{2}$ & $\pm\sqrt{3}$ \\
        \hline
            $4$ & $\pm\frac{\sqrt{2}}{2}$ & $\pm\frac{\sqrt{2}}{2}$ & $\pm1$ \\
        \hline
            $6$ & $\pm\frac{\sqrt{3}}{2}$ & $\pm\half$ & $\pm\frac{\sqrt{3}}{3}$ \\
        \hline
    \end{tabular}

    \vspace{1ex}
    \footnotesize \bf Figure 1: \normalfont All values of trigonometric functions \\ which are $\nth$ roots of rational numbers.

    \vspace{0.25ex}
    $b$ represents the denominator of the angle $\frac{a\pi}{b}$,
    
    \vspace{-0.8ex}
    and $a$ is any integer for which $\gcd (a,b) = 1$.
    \normalsize

    \vspace{6ex}

    {\includegraphics[width=12cm]{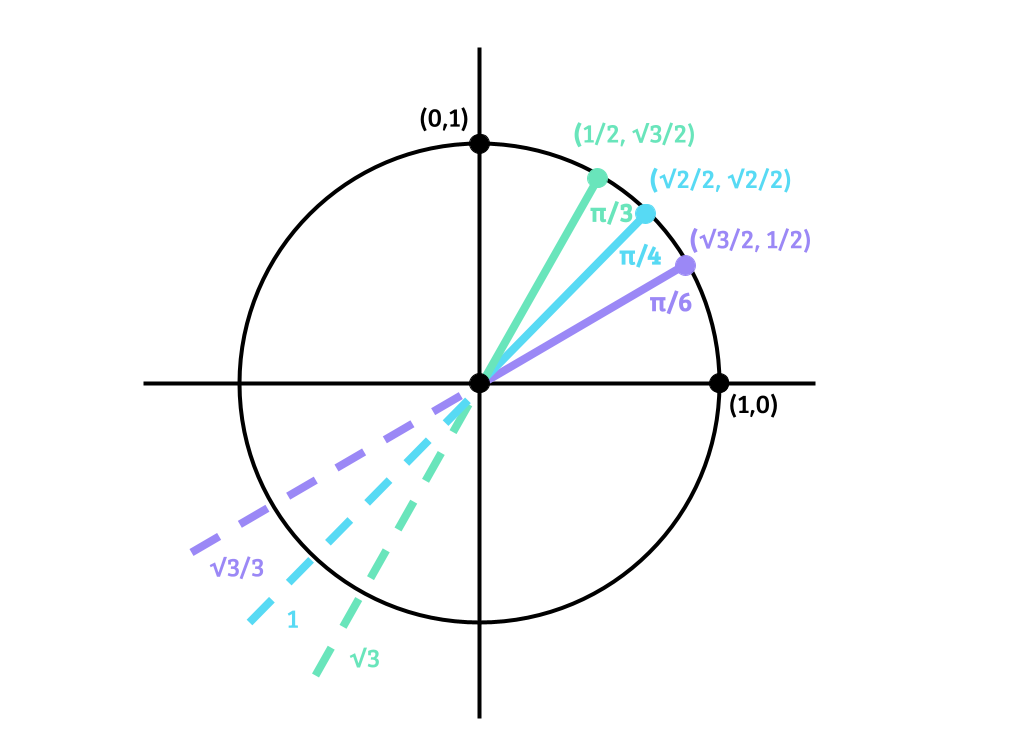}}

    \vspace{1ex}
    \footnotesize \bf Figure 2: \normalfont Some values of the unit circle.

    The coordinates to the top right \\ indicate cosine and sine pairs.

    The extended dashed lines indicate \\ the slope and are labelled by the \\ tangent values to the bottom left.
    
    \normalsize
\end{center}

\section{Background}
\subsection{Niven's Theorem and Existing Generalisations}

{
\begin{thm}[Niven's Theorem, \cite{niven}]\label{Niven}
    Suppose that both $r$ and $\text{cos}(r\pi) $ are elements of $ \mathbb Q$. Then $$\cos(r\pi) \in \left\{0, \pm\frac{1}{2}, \pm 1 \right\}.$$   
\end{thm}

\begin{proof}
    As $e^{r\pi i}$ is an algebraic integer for all $r\in\QQ$, we have that $$e^{r i\pi}+e^{-r i\pi} = 2\cos(r\pi)$$ is an algebraic integer; hence, as $\cos(r\pi)\in\QQ$, we must have $2\cos(r\pi) \in \mathbb Z$ if it is rational.
    
    Then, as $ 2\text{cos}\theta \in [-2,2]$, we must have that $2\text{cos}(\pi r) \in \{0, \pm 1, \pm 2\} $, and hence \[\cos(r\pi) \in \left\{0, \pm\frac{1}{2}, \pm 1 \right\}.\]
\end{proof}

\begin{rem}\label{rem:cos-to-sin}
    As we have the relation $\cos(r\pi)=\sin\left(r\pi + \frac{\pi}{2}\right)$, Niven's theorem also classifies values of $r$ such that $\sin(r\pi)$ is rational. Similarly, the question of finding rational $r$ such that $\sin(r\pi)$ has some algebraic property is equivalent to the same question for cosine.
\end{rem}

There have been several generalisations of Niven's Theorem, including the two following theorems, which are particularly relevant to this paper.

\begin{thm}[\cite{jahnel}]\label{quadratic}
    Among the quadratic irrationalities, only $\pm \frac{\sqrt{2}}{2}, \pm \frac{\sqrt{3}}{2}$ and $ \frac{\pm 1 \pm \sqrt{5}}{4}$ may be values of cos at rational angles. 
\end{thm}

Due to our focus on values of $\cos(r\pi)$ satisfying $\cos^n(r\pi) \in\QQ$, we are particularly interested in the values $\pm\frac{\sqrt{3}}{2}$, $\pm\frac{\sqrt{2}}{2}$, which satisfy $\cos^2(r\pi) \in\QQ$. However, the minimal polynomial of $\frac{\pm1 \pm\sqrt{5}}{4}$ is $4x^2\pm2x-1$ and so this value does not satisfy $\cos^n(r\pi) \in\QQ$.

\begin{thm}[\cite{nunn}]\label{tangent}
    Let $r\in\QQ$ and suppose $\tan(r\pi) \in \QQ$. Then $\tan(r\pi) \in \{0, \pm 1\}$.
\end{thm}

\subsection{Simple radicals and minimal polynomials}
\begin{defn}
    We call an algebraic number $\alpha$ a simple $n$\textsuperscript{th} radical if $\alpha^n\in\QQ$, and $\alpha^m\not\in\QQ$ for all $1\leq m<n$. We call $\alpha$ a simple radical if it is a simple $n$\textsuperscript{th} radical for some integer $n>0$. 
\end{defn}

Simple radicals have very simple minimal polynomials over $\QQ$, at least real simple radicals.

\begin{lem}\label{lem:minpolymn}
    Suppose that $\alpha \in\RR$ is a simple $\nth$ radical, with $\alpha^n=m\in\QQ$. The minimal polynomial of $\alpha$ over $\QQ$ is $x^n-m$.
\end{lem}

\begin{proof}
    Following \cite{koley,vahlen}, we have that $x^n - m \in \QQ[x]$ is reducible over $\QQ$ if and only if either \[m = k^t \, \text{ for some }\, t|n,\; t > 1,\; k\in \QQ\] or \[4|n \; \text{ and } \; m = -4k^4 \; \text{ for some }\; k \in \QQ.\]
    
    Based on our stipulation that $\alpha^N \notin \QQ$ for $N<n$, $N\in\NN$, we know that the first criterion does not hold. 
    
    Under the second criterion, we have $\alpha^{4t} = -4k^4$, or equivalently, $\alpha = \sqrt[4t]{-4k^4}$. This yields only imaginary values of $\alpha$, and is not applicable for our $\alpha\in\RR$. 
    
    Hence $x^n-m$ is irreducible. Therefore, for a simple $\nth$ radical $\alpha$, the minimal polynomial of $\alpha$ is $x^n-m$.
\end{proof}

Finally, we include a trivial property of minimal polynomials of rescaled algebraic numbers.

\begin{lem}\label{lem:minpoly-scale}
    If $\alpha$ is algebraic over $\QQ$ with monic minimal polynomial
    \[x^n + a_1x^{n-1} + \cdots + a_1 x+ a_0\]
    then, for any $r\in\QQ$, $r\alpha$ is algebraic with monic minimal polynomial
    \[x^n + r^{-1}a_1x^{n-1} + \cdots + r^{1-n}a_1 x + r^{-n}a_0.\]
\end{lem}

\section{Radical values of sine and cosine}
As noted earlier in Remark \ref{rem:cos-to-sin}, the problem of finding rational $r$ such that $\sin(r\pi)$ is a simple radical is equivalent to that for $\cos(r\pi)$. As such, we will consider only the cosine case. To do so, we consider properties of the minimal polynomial of cosine, presented by Berger \cite{berger} as derived from Lehmer \cite{lehmer} and Kunz \cite{kunz}.

\begin{prop}[Proposition 2.7, \cite{berger}]
     For every integer $n \geq 0$ there exists a unique monic polynomial $R_n(x)\in \ZZ[x]$ such that $$x^n + x^{-n} = R_n\left(x + x^{-1}\right), \; \forall x \in \CC \backslash \{0\}.$$
\end{prop}

Using these $R_n(x)$, we can construct the minimal polynomial of $\cos(r\pi)$ as follows. Given $n\geq 0$, we write the $n$\textsuperscript{th} cyclotomic polynomial as
$$\Phi_n(x) = \sum_{j=0}^{\varphi(n)} a(j,n)x^{\varphi(n)-j}.$$ 
For convenience, we extend $a(j,n)$ to be fully defined for all $n\in\NN, j\geq0$, by setting $a(j,n)=0$ whenever $j\geq\varphi(n)$.

We define a polynomial $P_n(x)\in \ZZ[x]$ by $P_1(x)=x+2$ and
\[P_n(x) = \sum_{j=0}^{\half\varphi(2n)-1}a(j,2n)R_{\half\varphi(2n)-1}(x)+a\left(\half\varphi(2n),2n\right) \text{\quad} \forall n\geq2.\]

\begin{prop}[Lemma 2.8, \cite{berger}]\label{prop:min_poly_properties} 
    Suppose $r=\frac{a}{b}\in\QQ$, with $a,\, b$ coprime integers, $b>0$. Then $P_b(x)$ is the monic minimal polynomial of $2(-1)^{1+a}\cos(r\pi)$ over $\QQ$. Furthermore, the degree of $P_n(x)$ is

\[\begin{aligned}
  p_n: =\deg P_n(x)  = &\begin{cases}
    1 & \text{if } n = 1 ,\\
    \frac{1}{2} \varphi(2n) & \text{if } n \geq 2,
  \end{cases} \\
  =&\begin{cases}
    1 & \text{if } n = 1, \\
    \varphi(n) & \text{if } n \geq 2 \text{ even}, \\
    \frac{1}{2} \varphi(n) & \text{if } n \geq 2 \text{ odd}.
  \end{cases}
\end{aligned}
\]

The constant term $P_n(0)$ is given by
 \[\text{\quad\quad\quad\quad\quad\quad\quad\quad\quad\quad}\;\; |P_n(0)| = \begin{cases}
    2 & \tif n = 1, \\
    0 & \tif n = 2, \\    
    p & \tif n = 2p^j \text{ for some $p$ prime and $j \in \NN$},\\    
    1 & \text{ otherwise}.
\end{cases}\]
\end{prop}

Note that, essentially by definition, we can relate $P_n(x)$ very explicitly to the cyclotomic polynomials.

\begin{prop}\label{prop:min_poly_to_cyclotomic}
    For all $n\geq 2$,
    \[\Phi_{2n}(x) = x^{p_n}P_n(x+x^{-1}).\]
\end{prop}

\subsection{Simple radical values of cosine}
If $\cos(r\pi)$ is a simple $n$\textsuperscript{th} radical, Lemma \ref{lem:minpolymn} tells us that its minimal polynomial is $x^n-m$ for some $m\in\QQ$. Proposition \ref{prop:min_poly_properties} tells us exactly the minimal polynomial of $\cos(r\pi)$, and so $\cos(r\pi)$ is a simple $n$\textsuperscript{th} radical if and only if
\[P_b(x) = x^n - \left(2(-1)^{1+a}\right)^nm\]

We first note that if $b\neq 2p^j$ for some prime $p$, then we are in one of three cases. If $b=1$, then $\cos(r\pi)=\pm 1$. If $b=2$, then $\cos(r\pi)=0$. Finally, if $b\neq 1,2,2p^j$ and is a simple radical, we must have that
\[\cos^n(r\pi) = \pm\left(\half\right)^n\]
so $\cos(r\pi)=\pm \half$. These are precisely the rational values arising in Niven's Theorem. As such, we can only get interesting radical values if $b=2p^j$, and so we will assume $b$ is of this form from here onwards.

Comparing the degrees of these polynomials, we can immediately conclude the following.
\begin{lem}
    Suppose $r=\frac{a}{b}$ where $b=2p^j$, and suppose that $\cos(r\pi)$ is a simple $n$\textsuperscript{th} radical. Then $n=2^j$ if $p=2$, and $n=p^{j-1}(p-1)$ if $p$ is an odd prime.
\end{lem}

Similarly, comparing the constant coefficients, we obtain constraints on the constant term $m$.
\begin{cor}\label{cor:nm}
    Let $r=\frac{a}{b}$ and $b=2p^j$, and suppose $\cos(r\pi)$ is a simple $\nth$ radical. The constant $m = \cos^n(r\pi)$ is $m=\pm2^{1-2^j}$ when $p=2$, and $m=\pm {p}\cdot {2^{p^{j-1}(1-p)}}$ when $p$ is an odd prime.
\end{cor}

\begin{proof}
    From Proposition \ref{prop:min_poly_properties}, we know that $|m| = \frac{p}{2^{n}}$ when $b=2p^j$.

    First consider $p=2$, so that $n=2^j$. Then $$|m|=\frac{2}{2^{2^j}}=2\cdot2^{-2^j}=2^{1-2^j}.$$

    Now let $p$ be some odd prime, so that $n=p^{j-1}(p-1)$. Then
    
    $$|m|=\frac{p}{2^{p^{j-1}(p-1)}} = p\cdot 2^{-p^{j-1}(p-1)} = p\cdot 2^{p^{j-1}(1-p)}.$$
\end{proof}

Proposition \ref{prop:min_poly_to_cyclotomic} then gives us an explicit formula for certain cyclotomic polynomials if $\cos(r\pi)$ is a simple radical. Let us first consider the case where $b=2p^j$ for an odd prime $p$. Then, if $\cos(r\pi)$ is a simple radical, we must have that
$$\Phi_{4p^j}(x) = x^{p^{j-1}(p-1)}\left(\left(x+x^{-1}\right)^{p^{j-1}(p-1)}\pm p\right).$$ 
As we are assuming this has a real root, we must in fact have
$$\Phi_{4p^j}(x) = x^{p^{j-1}(p-1)}\left(\left(x+x^{-1}\right)^{p^{j-1}(p-1)}- p\right).$$ 
In particular, we must have $\Phi_{4p^j}(1) =2^{p^{j-1}(p-1)} - p$. We proceed to find an explicit polynomial expression of the $4p^j$-th cyclotomic polynomial.

\begin{lem}\label{phi4pj}
    For all odd primes $p$, $\Phi_{4p^j}(x)= \frac{x^{2p^j}+1}{x^{2p^{j-1}}+1}$.
\end{lem}
\begin{proof}
    Recall that $\prod_{d|n}\Phi_d(x) = x^n-1$. Hence
    \begin{align*}
        x^{2p^j}-1 &=\prod_{d|2p^j}\Phi_d(x) \\
        &=\prod_{d|p^j}\Phi_d(x)\Phi_{2d}(x) \\
        &=\left(x^{p^j}-1\right)\prod_{d|p^j}\Phi_{2d}(x).
    \end{align*}
    As
    $$x^{2p^j}-1 = \left(x^{p^j}-1\right)\left(x^{p^j}+1\right),$$ we have $$x^{p^j}+1=\prod_{d|p^j}\Phi_{2d}(x).$$    
    Similarly,
    \begin{align*}
        x^{4p^j}-1&=\prod_{d|4p^j}\Phi_d(x) \\
        &=\prod_{d|p^j}\Phi_d(x)\Phi_{2d}(x)\Phi_{4d}(x) \\
        &=\left(x^{p^j}-1\right)\left(x^{p^j}+1\right) \prod_{d|p^j} \Phi_{4d}(x).
    \end{align*}
    and 
    $$x^{4p^j}-1 = \left(x^{2p^j}-1\right)\left(x^{2p^j}+1\right) = \left(x^{p^j}-1\right)\left(x^{p^j}+1\right) \left(x^{2p^j}+1\right),$$
    so we must have that
    $$x^{2p^j}+1=\prod_{d|p^j}\Phi_{4d}(x).$$
    Thus
    \[\Phi_{4p^j}(x) = \frac{x^{2p^j}+1}{\prod_{d|p^{j-1}}\Phi_{4d}(x)} = \frac{x^{2p^j}+1}{x^{2p^{j-1}}+1}.\]
\end{proof}

We therefore have that
\[2^{p^{j-1}(p-1)}-p = \Phi_{4p^j}(1) = \frac{1^{2p^j}+1}{1^{2p^{j-1}}+1} = \frac{2}{2} = 1\]
if $\cos(r\pi)$ is a simple radical.

\begin{prop}\label{oddprimes}
    The only odd prime $p$ and positive integer $j$ satisfying $$2^{p^{j-1}(p-1)} - p = 1$$ are $p=3,j=1$.
\end{prop}

\begin{proof}
    Let $p$ be some odd prime and let $j\in\NN$ be a positive integer.
    
   If $p=3$, we can explicitly solve\begin{align*}
        1 &= 2^{3^{j-1}(3-1)} - 3 
        \end{align*}
    to find $j=1$.

    Now suppose $p\geq 5$. Then it is not difficult to show by induction that
    \[2^{p^{j-1}(p-1)}-p \geq 2^{p-1}-p > (p+1)-p = 1\]
    and so we cannot have any other solutions to the given equation.
\end{proof}

\begin{cor}\label{cor:odd_primes}
    Suppose $r=\frac{a}{b}$ for coprime $a,b$, with $b=2p^j$ for some odd prime $p$. Then if $\cos(r\pi)$ is a simple radical, $p=3$, $j=1$, and $\cos(r\pi)=\pm\frac{\sqrt{3}}{2}$.
\end{cor}

Now we consider the case of $p=2$ via a similar approach. As before, if $\cos(r\pi)$ is a simple radical, we must have that
\[\Phi_{2b}(x)=\Phi_{2(2^{j+1})}(x) = x^{2^j}\left(\left(x+x^{-1}\right)^{2^j}-2\right)\]
and so
\[\Phi_{2^{j+2}}(1)=2^{2^{j}}-2\]
A very similar calculation gives us the following evaluation of $\Phi_{2^{j+2}}(x)$.
\begin{lem}
    For all $j\geq 0$, $\Phi_{2^{j+2}}(x)=x^{2^{j+1}}+1$, and $\Phi_{2^{j+2}}(1)=2$.
\end{lem}

\begin{cor}\label{cor:even_primes}
    Suppose $r=\frac{a}{b}$ for coprime $a,b$, with $b=2^{j+1}$. Then if $\cos(r\pi)$ is a simple radical, $j=1$, and $\cos(r\pi)=\pm\frac{\sqrt{2}}{2}$.
\end{cor}
\begin{proof}
    If these assumptions hold, we must have that
    \[2^{2^j}-2 = 2\]
    which we can directly solve to find $j=1$, implying $\cos(r\pi)=\cos\left(\frac{a\pi}{4}\right)=\pm\frac{\sqrt{2}}{2}$.
\end{proof}

Combining the results of this section, we therefore have
\begin{thm}\label{cosineCase}
    Suppose that $\cos^n(r\pi)$ is rational for $r\in\QQ$, $n\in\NN$. Then \[\cos(r\pi) \in \left\{0, \pm 1,\pm\half,\pm\frac{\sqrt{3}}{2},\pm\frac{\sqrt{2}}{2}\right\}.\]
\end{thm}

\begin{cor}
    Suppose that $\sin^n(r\pi)$ is rational for $r\in\QQ$, $n\in\NN$. Then \[\sin(r\pi) \in \left\{0, \pm 1,\pm\half,\pm\frac{\sqrt{3}}{2},\pm\frac{\sqrt{2}}{2}\right\}.\]
\end{cor}

Notably, we find that if $\cos(r\pi)$ is an $\nth$ simple radical, then $n\leq 2$. A more conceptual explanation for this is given in the final section.

\section{Radical values of tangents}\label{sec:tangent}
While Galieva \& Galyautdinov provide a method of calculating these polynomials recursively \cite{russian}, there is no simple description known to the authors, and so the same arguments cannot be used in the case of tangent. However, using standard properties of degrees of field extensions and norms of algebraic numbers, we can also determine all possible radical values of $\tan(r\pi)$ for $r\in\QQ$.

Suppose that $\tan(r\pi)$ is a simple $n$\textsuperscript{th} radical, with $\tan^n(r\pi)=m\in\QQ$. This is equivalent to
\[\sin(r\pi) = \cos(r\pi)\sqrt[n]{m}\]
for an appropriate (real) $n$\textsuperscript{th} root of $m$.

Define a field extension $L=\QQ[\sin(r\pi),\cos(r\pi)]$ of $\QQ$, which must also contain $\sqrt[n]{m}$. Taking the norm of the above equation, we must have
$$N^L_\QQ\left(\sin(r\pi)\right) =  N^L_\QQ\left(\cos(r\pi)\sqrt[n]{m}\right) =  N^L_\QQ\left(\cos(r\pi)\right)N^L_\QQ\left(\sqrt[n]{m}\right).$$

Let $S(x)$ and $C(x)$ be the minimal polynomials of $\sin(r\pi)$ and $\cos(r\pi)$ over $\QQ$ respectively, and let $c_S$, $c_C$ denote their constant terms. The minimal polynomial of $\sqm$ over $\QQ$ is $x^n-m$, as shown in Lemma \ref{lem:minpolymn}.

Next, let $L_S=\QQ\left(\sin(r\pi)\right)$, $L_C=\QQ\left(\cos(r\pi)\right)$, and $L_m = \QQ\left(\sqm\right)$ be extensions of $\QQ$ generated by $\sin(r\pi)$, $\cos(r\pi)$, and $\sqm$ respectively. The equality of norms then implies the following.

\begin{lem}\label{constants}
    $\pm c_S^{[L:L_S]} = c_c^{[L:L_C]}m^{[L:L_m]}$.
\end{lem}
\begin{proof}
    Recall that, for $\alpha\in L$,
    $$N^L_\QQ(\alpha) = \left(\prod_{i=1}^k\sigma_i(\alpha)\right)^{[L:\QQ(\alpha)]},$$
    where we take the product over all possible embeddings $\sigma_i:L\to \CC$.
    We also have that $\prod_{i=1}^k\sigma_i(\alpha)$ is the constant term of the minimal polynomial of $\alpha$ over $\QQ$, up to a sign. The result then follows.
\end{proof}

From the identity
\[\sin^2(r\pi)+\cos^2(r\pi)=1\]
we can conclude that 
\[[L:L_S]\leq 2\quad\text{and}\quad [L:L_C]\leq 2.\]

From Proposition \ref{prop:min_poly_properties}, we must have that
\[[L_C:\QQ] = p_b = \begin{cases}
    1 & \tif b = 1, \\
    \varphi (b) & \tif b \geq 2 \text{ even}, \\    
    \half \varphi (b) & \tif b \geq 2 \text{ odd}.
\end{cases} \]

Similarly, from the minimal polynomial of $\sin(rx)$ \cite{berger}[Lemma 2.10], we must have that
\[[L_S:\QQ] = \begin{cases}
    1 & \tif b = 2, \\
    \half \varphi (b) & \tif b\in 2+4\ZZ \text{ and } b \geq 6, \\
    \varphi (b) & \text{ otherwise}.
\end{cases} \]

\begin{prop}\label{prop:tan_equal_case}
    Suppose $r=\frac{a}{b}$ for coprime integers $a,b$, with $b>0$, and that $\tan(r\pi)$ is a simple $n$\textsuperscript{th} radical. Let $L,\, L_S,\, L_C$ be as above, and suppose that $[L:L_S]=[L:L_C]$. Then $b\in\{1,4\}$ and 
    \[\tan(r\pi)\in\{0,\pm 1\}.\]
\end{prop}
\begin{proof}
    From tower law, $[L:L_S]=[L:L_C]=d$ if and only if 
    \[[L_S:\QQ]=[L_C:\QQ].\]
    Comparing the degrees of these extensions across possible $b$, we see that this equality implies that $b\in \{1,2\}$ or $b\in 4\ZZ$. If $b=1$, then $\tan(r\pi)=0$. If $b=2$, then $\tan(r\pi)$ is undefined and can be ignored. When $b\in 4\ZZ$, $P_b(x)$ is the minimal polynomial of both one of $\{ \cos(r\pi), -\cos(r\pi)\}$ and one of $\{ \sin(r\pi),-\sin(r\pi)\}$ \cite{berger}. Hence $c=c_S=\pm c_C$, and so Lemma \ref{constants} implies that
    \[c^d=\pm c^d m^{[L:L_m]}.\]
    If $b\in 4\ZZ$, then $c\neq 0$ and so we must have $m=\pm 1$, which implies we must in fact have $b=4$, from which the claim follows.
\end{proof}

\begin{prop}\label{prop:tan_unequal_case}
    Suppose $r=\frac{a}{b}$ for coprime integers $a,b$, with $b>0$, and that $\tan(r\pi)$ is a simple $\nth$ radical. Let $L,\, L_S,\, L_C$ be as above, and suppose that $[L:L_S]\neq[L:L_C]$. Then $b\in\{3,6\}$ and 
    \[\tan(r\pi)\in\{\pm\sqrt{3}^{\pm 1}\}\]
\end{prop}
\begin{proof}
    We have two cases: $[L:L_S]=1$, $[L:L_C]=2$ and $[L:L_S]=2$, $[L:L_C]=1$. We will only handle the first case, as the second is similar.

    If $[L:L_S]=1$, then $\cos(r\pi)\in L=L_S$, and so there exists a polynomial $f(x)\in\QQ[x]$ such that
    \[\cos(r\pi)=f(\sin(r\pi)).\]
    Furthermore $\deg f\geq 2$, as $L\neq L_C$. As $(\cos(r\pi)^2=1-(\sin(r\pi))^2$, we have that $\sin(r\pi)$ is a root of
    \[(f(x))^2 +x^2 -1.\]
    Recall that $\sin(r\pi)$ has minimal polynomial $S(x)$. By minimality, we must have that $2\deg f>\deg S$, and so we can find unique $g(x),h(x)\in \QQ[x]$ with $\deg(h)<\deg(S)$ such that 
    \[f^2(x) = S(x)g(x) + h(x) \]
    In particular
    \[f^2(\sin(r\pi)) = S(\sin(r\pi))g(\sin(r\pi)) + h(\sin(r\pi)) = h(\sin(r\pi))\]
    and so
    \[ 1- \text{sin}^2(r\pi) = h(\sin(r\pi)).\]
Thus $\sin(r\pi)$ is a root of $h(x) + x^2 -1$. If $\deg S>2$, then this contradicts the minimality of $S(x)$. Tower law tells us that 
\[\deg S = [L_S:\QQ] = [L:\QQ] = 2[L_C:\QQ] = 2\deg C \geq 2\]
and hence we must have $\deg S=2$, $\deg C=1$.

Similarly, if $[L:L_S]=2$, $[L:L_C]=1$, then $\deg S=1$ and $\deg C=2$. Thus, if $\tan(r\pi)$ is a simple radical, one of $\sin(r\pi)$ and $\cos(r\pi)$ is rational and the other is quadratic. As such, $\tan(r\pi)$ must be quadratic. Running through all possible $b$ such that one of $\sin(r\pi)$ and $\cos(r\pi)$ is quadratic and one is rational, we find that $(\sin(r\pi),\cos(r\pi))$ is one of the pairs
\[\left(\pm \frac{\sqrt{3}}{2}, \pm \half\right),\quad \text{or}\quad \left(\pm \half,\pm \frac{\sqrt{3}}{2}\right)\]
from which the claim follows.
\end{proof}

Combining all the results of this section, we find

\begin{thm}
    Suppose that $\tan^n(r\pi)$ is rational for $r\in\QQ$, $n\in\NN$. Then \[\tan(r\pi) \in \left\{0, \pm 1,\pm\sqrt{3}^{\pm 1}\right\}.\]
\end{thm}

\section{A Galois theoretic explanation}
While the primary goal of this article was to discuss radical values of trigonometric functions using as close to elementary methods as possible, it would be unfulfilling not to address a more direct and conceptual route to our results using notions from Galois theory. This relies on two key results.

Recall first that an extension $K/\QQ$ is called normal if every irreducible polynomial $f(x)\in\QQ[x]$ with a root in $K$ splits as a product of linear factors in $K[x]$.

\begin{lem}\label{lem:normal_extensions}
    Let $x^n-m\in \QQ[x]$ have a real root $\alpha\in\RR$. Then $\QQ[\alpha]/\QQ$ is not normal if $n>2$.
\end{lem}
\begin{proof}
    Every such polynomial has at least one complex root given by $\alpha e^{\frac{2\pi}{n}}$, which cannot be in the real field $\QQ[\alpha]\subset\RR$.
\end{proof}

Next we need a more explicit application of the Galois correspondence

\begin{lem}\label{lem:subfields_of_abelian}
  Let $K/\QQ$ be a Galois extension, with abelian Galois group $G$. Then every intermediate field $\QQ\subset F\subset K$ is a normal extension of $\QQ$.  
\end{lem}
\begin{proof}
    The Galois correspondence tells us that every intermediate field $F$ is in bijection is a subgroup $H^F$ of $G$, and furthermore that $F/\QQ$ is a normal extension if and only if $H^F$ is a normal subgroup of $G$. Every subgroup of an abelian group is normal, and so the result follows.
\end{proof}

\begin{prop}
    Suppose that $\cos(r\pi)$ is a simple $n$\textsuperscript{th} radical. Then $n\leq 2$.
\end{prop}
\begin{proof}
    Note that $\cos(r\pi)=\half\left(e^{r\pi}+e^{-r\pi}\right)$, and so
    \[\QQ\subset \QQ[\cos(r\pi)]\subset \QQ[e^{r\pi}].\]
    It is a standard result that $\Gal(\QQ[e^{r\pi}]/\QQ)$ is a cyclic group, and in particular abelian. Thus, Lemma \ref{lem:subfields_of_abelian} implies that $\QQ[\cos(r\pi)]/\QQ$ is a normal extension. Lemma \ref{lem:normal_extensions} then implies that if $\cos(r\pi)$ is a simple $n$\textsuperscript{th} radical, $n\leq 2$.
\end{proof}

Theorem \ref{cosineCase} then follows immediate from Theorem \ref{quadratic}. Similarly, as $\QQ[\sin(r\pi)]$ and $\QQ[\tan(r\pi)]$ are subfields of the cyclotomic field $\QQ[e^{r\pi}]$, $\sin(r\pi)$ and $\tan(r\pi)$ can be simple $n$\textsuperscript{th} radicals only for $n\leq 2$, from which we can make the same conclusions.

\printbibliography
\addcontentsline{toc}{section}{References}
\Addresses
\end{document}